\documentclass[12pt,twoside,leqno]{article}

\usepackage{amssymb}
\usepackage{amsmath}
\usepackage{bbm}
\usepackage{mathrsfs}

\makeatletter
\gdef\th@mychange{\normalfont\slshape
   \def\@begintheorem##1##2{\item
        [\hskip\labelsep \theorem@headerfont ##2. ##1  \,--\!--\!--\!--  ]}%
 \def\@opargbegintheorem##1##2##3{%
   \item[\hskip\labelsep \theorem@headerfont ##2. ##1\ {\upshape(}##3{\upshape)}. \,-----  ]}}
\makeatother

\RequirePackage{theorem}

\theoremstyle{mychange}

{\theorembodyfont{\rmfamily}\newtheorem{ttt}{\!\!\!}[section]}
{\theorembodyfont{\rmfamily}\newtheorem{ttt-s}{\!\!\!}[subsection]}

{\theorembodyfont{\rmfamily}\newtheorem{rem}[ttt]{Remark.}}
{\theorembodyfont{\rmfamily}\newtheorem{remso}[ttt]{Remarks}}
{\theorembodyfont{\rmfamily}\newtheorem{rems}[ttt]{Remarks.}}
{\theorembodyfont{\rmfamily}}

{\theorembodyfont{\itshape}}

{\theorembodyfont{\rmfamily}\newtheorem{c81-s}[ttt-s]{Congruences modulo 81.}}
{\theorembodyfont{\rmfamily}\newtheorem{sit-s}[ttt-s]{The Situation.}}
{\theorembodyfont{\rmfamily}\newtheorem{prog-s}[ttt-s]{Our Program.}}
{\theorembodyfont{\rmfamily}\newtheorem{prob-s}[ttt-s]{Programmer's position.}}
{\theorembodyfont{\rmfamily}\newtheorem{sol-s}[ttt-s]{A way out.}}
{\theorembodyfont{\rmfamily}\newtheorem{why-s}[ttt-s]{Why is this algorithm faster?}}

{\theorembodyfont{\rmfamily}\newtheorem{remo-s}[ttt-s]{Remark}}
{\theorembodyfont{\rmfamily}\newtheorem{rem-s}[ttt-s]{Remark.}}
{\theorembodyfont{\rmfamily}\newtheorem{rems-s}[ttt-s]{Remarks.}}
{\theorembodyfont{\itshape}\newtheorem{theo-s}[ttt-s]{Theorem.}}

\usepackage{fancyhdr}
\newcounter{ABC}
\newenvironment{ABC}{\begin{list}{\rm \Alph{ABC}. }{\usecounter{ABC} \leftmargin=0.0pt \labelsep=0.0pt \listparindent=0.0pt \labelwidth=0.0pt \parsep=\smallskipamount \itemsep=0.0pt \topsep=0.0pt \partopsep=\smallskipamount}}{\end{list}}

\newcounter{III}
\newenvironment{III}{\begin{list}{\rm \Roman{III}. }{\usecounter{III} \leftmargin=0.0pt \labelsep=0.0pt \listparindent=0.0pt \labelwidth=0.0pt \parsep=\smallskipamount \itemsep=0.0pt \topsep=0.0pt \partopsep=\smallskipamount}}{\end{list}}

\newcounter{abc}
\newenvironment{abc}{\begin{list}{\rm \alph{abc}) }{\usecounter{abc} \leftmargin=0.0pt \labelsep=0.0pt \listparindent=0.0pt \labelwidth=0.0pt \parsep=\smallskipamount \itemsep=0.0pt \topsep=0.0pt \partopsep=\smallskipamount}}{\end{list}}

\newcounter{iii}
\newenvironment{iii}{\begin{list}{\rm \roman{iii}) }{\usecounter{iii} \leftmargin=0.0pt \labelsep=0.0pt \listparindent=0.0pt \labelwidth=0.0pt \parsep=\smallskipamount \itemsep=0.0pt \topsep=0.0pt \partopsep=\smallskipamount}}{\end{list}}

\newcommand{\bbF}{{\mathbbm F}}
\newcommand{\bbN}{{\mathbbm N}}
\newcommand{\bbP}{{\mathbbm P}}
\newcommand{\bbZ}{{\mathbbm Z}}

\newcommand{\notd}{\mathord{\nmid}}

\renewcommand{\mod}{\text{ {\rm mod} }}

\def\pmod#1{\nobreak\ifinner\mkern8mu\else\mkern18mu\fi (\text{\rmfamily\upshape mod}\,\,#1)}

\renewcommand{\thefootnote}{\arabic{footnote}}
\author{Andreas-Stephan Elsenhans${}^{1,2}$ and J\"org Jahnel${}^2$}
\date{}

\title{The Diophantine Equation
$x^4 + 2 y^4 = z^4 + 4 w^4$--- \\
a number of improvements}

\begin{document}

\maketitle

\begin{abstract}
The quadruple
$(1\,484\,801, 1\,203\,120, 1\,169\,407, 1\,157\,520)$
already known is essentially the only non-trivial solution of the Diophantine equation 
$x^4 + 2 y^4 = z^4 + 4 w^4$
for
$|x|$,
$|y|$,
$|z|$,
and
$|w|$
up to one hundred~million. We~describe the algorithm we used in order to establish this result, thereby explaining a number of improvements to our original approach~\cite{Joerg_und_Stephan}.
\end{abstract}

\section{Introduction}
\footnotetext[1]{The first author was partially supported by a Doctoral Fellowship of the Deutsche Forschungsgemeinschaft~(DFG).}
\footnotetext[2]{The computer part of this work was executed on the Linux~PCs of the Gau \ss\ Laboratory for 
Scientific Computing at the G\"ottingen Mathematical Institute. 
Both authors are grateful to Prof.~Y.~Tschinkel for the permission to use these 
machines as well as to the system administrators for their~support.}

\thispagestyle{empty}

\begin{ttt}
In \cite{Joerg_und_Stephan}, we described a systematic method to search efficiently for all solutions of a Diophantine equation of the form 
$$f (x_1, \ldots, x_n) = g (y_1, \ldots, y_m)$$
which are contained within the
$(n + m)$-dimensional
cube 
$$\{ (x_1, \ldots, x_n, y_1, \ldots, y_m) \in \bbZ^{n+m} \mid |x_i|, |y_i| \leq B \}.$$
The expected running-time of this algorithm
is~$O(B^{\max \{ n, m \}})$.
\end{ttt}

\begin{ttt}
The basic idea is as~follows.\medskip

{\bf Algorithm H.}
\begin{iii}
\item
Evaluate
$f$
on all points of the
$n$-dimensional
cube
$\{ (x_1, \ldots, x_n) \in \bbZ^n \mid |x_i| \leq B \}$.
Store the values within a
set~$L$.
\item
Evaluate
$g$
on all points of the cube
$\{ (y_1, \ldots, y_m) \in \bbZ^m \mid |y_i| \leq B \}$
of
dimension~$m$.
For~each value start a search in order to find out whether it occurs
in~$L$.
When~a coincidence is detected, reconstruct the corresponding values of
$x_1, \ldots, x_n$
and output the~solution.
\end{iii}
\end{ttt}

\begin{rems} ~
\begin{abc}
\item
In fact, we are interested in the very particular Diophantine equation
$x^4 + 2 y^4 = z^4 + 4 w^4$
which was suggested by Sir Peter Swinnerton-Dyer. It~is unknown whether this equation admits finitely or infinitely many primitive~solutions. If~their number were actually finite then this would settle a famous open problem in the arithmetic of
$K3$~surfaces~\cite[Problem/Question~6.a)]{PT}.
\addtocounter{abc}{1}
\begin{iii}
\item[b.i) ]\addtocounter{iii}{1}
In the form stated above, the main disadvantage of Algorithm~H is that it requires an enormous amount of~memory. Actually,~the set 
$L$
is too big to be stored in the main memory even of our biggest computers, already when the value
of~$B$
is only moderately~large.

For that reason, we introduced the idea of {\em paging}. We choose a {\em page~prime}
$p_p$
and work with the sets
$L_r := \{ s \in L \mid s \equiv r \pmod {p_p} \}$
for
$r = 0, \ldots, p_p - 1$,~separately.
At the cost of some more time spent on initializations, this yields a~reduction of the memory space required by a factor
of~$\smash{\frac1{p_p}}$.
\item
The sets
$L_r$
were implemented in the form of a hash~table with open~addressing.
\item
It is possible to achieve a further reduction of the running-time and the memory required by making use of some obvious congruence~conditions
modulo~$2$
and~$5$.
\end{iii}
\item
Precisely ten primitive solutions of the Diophantine equation
$x^4 + 2 y^4 = z^4 + 4 w^4$
are known up to~now. Among~them, there are the two obvious ones
$(\pm1 \!:\! 0 \!:\! \pm1 \!:\! 0)$.

Furthermore,~by an implementation of Algorithm~H, the non-obvious solutions
$(\pm1\,484\,801 \!:\! \pm1\,203\,120 \!:\! \pm1\,169\,407 \!:\! \pm1\,157\,520)$
were~found. We~searched through the hypercube
$\{ (x, y, z, w) \in \bbZ^4 \mid |x|, |y|, |z|, |w| \leq 2.5 \cdot 10^6 \}$.
Details are given in~\cite{Joerg_und_Stephan}.
\end{abc}
\end{rems}

\begin{ttt}
The goal of this note is to describe an improved implementation of Algorithm~H which we used in order to find all solutions of
$x^4 + 2 y^4 = z^4 + 4 w^4$
contained within the hypercube
$\{ (x, y, z, w) \in \bbZ^4 \mid |x|, |y|, |z|, |w| \leq 10^8 \}$.

Unfortunately, our result is not very~spectacular. There~is no new primitive~solution.
\end{ttt}

\section{More Congruences}

\begin{ttt-s}
The most obvious way to further reduce the size of the
sets~$L_r$
and to increase the speed of Algorithm~H is to find further congruence~conditions for solutions and evaluate
$f$
and~$g$
only on points satisfying these~conditions. As~the equation, we are interested in, is homogeneous, it is sufficient to restrict consideration to primitive~solutions.
\end{ttt-s}

\begin{ttt-s}
It should be noticed, however, that this idea is subject to strict~limitations. If we were using the most naive
$O (B^{n+m})$-algorithm
then, for more or less every
$l \in \bbN$,
the congruence
$f (x_1, \ldots, x_n) \equiv g (y_1, \ldots, y_m) \pmod l$
caused a reduction of the number of
$(n+m)$-tuples
to be~checked. For~Algorithm~H, however, the situation is by far less~fortunate.

One~may gain something only if there are residue classes
$(r \mod l)$
which are represented
by~$f$,
but not
by~$g$,
or vice~versa. Values,~the residue class of which is not represented
by~$g$,
do not need to be stored
into~$L_r$.
Values,~the residue class of which is not represented
by~$f$,
do not need to be searched~for.

Unfortunately,~if
$l$
is~prime and not very small then the Weil conjectures ensure that all residue classes
modulo~$l$
are~represented by both
$f$
and~$g$.
In~this case, the idea fails~completely. The~same is, however, not true for prime~powers~$l = p^k$.
Hensel's~Lemma does not work when all partial derivatives
$\smash{\frac{\partial{f}}{\partial{x_i}} (x_1, \ldots, x_n)}$,
respectively~$\smash{\frac{\partial{g}}{\partial{y_i}} (y_1, \ldots, y_m)}$,
are divisible
by~$p$.
This~makes it possible that certain residue classes
$(r \mod p^k)$
are not representable although
$(r \mod p)$~is.
\end{ttt-s}
 
\subsection{The prime~$5$. Congruences modulo~$625$}

\begin{ttt-s}
In~\cite{Joerg_und_Stephan}, we made use of the fact that
$y$
is always divisible
by~$5$.
However, at this point, one can do a~lot~better. When~one takes into consideration that
$a^4 \equiv 1 \pmod 5$
for
every~$a \in \bbZ$
not divisible
by~$5$,
a systematic inspection shows that there are actually two~cases.

Either,~$5 | w$.
Then,~$5 \notd x$
and~$5 \notd z$.
Or,~otherwise,
$5 | x$.
Then,~$5 \notd z$
and~$5 \notd w$.
Note~that, in the latter case, one indeed has 
$z^4 + 4 w^4 \equiv 1 + 4 \equiv 0 \pmod 5$.
\end{ttt-s}

\begin{ttt-s}
{\bf The Case $5 | w$.}
We call this case ``N'' and use the letter~N at a prominent position in the 
naming of the relevant files of source~code. N~stands for~``normal''. 
To~consider this case as the ordinary one is justified by the fact that 
all primitive solutions known actually belong to~it. Note,~however, that 
we have no theoretical reason to believe that this case should in whatever 
sense be better than the other~one.

In case~N, we rearrange the equation to
$f_N (x, z) = g_N (y, w)$
where
$$f_N (x, z) := x^4 - z^4 \qquad {\rm and} \qquad g_N (y, w) := 4 w^4 - 2 y^4.$$
As
$y$
and
$w$
are both divisible
by~$5$,
we get
$g_N (y, w) = 4 w^4 - 2 y^4 \equiv 0 \pmod {625}$.
Consequently,~$f_N (x, z) \equiv 0 \pmod {625}$.

This yields an enormous reduction of the
set~$L_r$.
To~see this, recall
$5 \notd x$
and~$5 \notd z$.
That~means,
for~$x$,
there are precisely
$\varphi(625)$
possibilities
in~$\bbZ/625\bbZ$.
Further,~for each such value, the congruence
$z^4 \equiv x^4 \pmod {625}$
may not have more than four~solutions. All~in all, there are
$4 \cdot \varphi(625) = 2\,000$
possible pairs
$(x, z) \in \left( \bbZ/625\bbZ \right)^2$.

Further, these pairs are very easy to find,~computationally. The~fourth roots of unity 
modulo~$625$
are
$\pm1$
and~$\pm182$.
For~each
$x \in \bbZ/625\bbZ^{\textstyle *}$,
put
$z := (\pm x \mod 625)$
and
$z := (\pm 182x \mod 625)$.

We store the values of
$f_N$
into the set
$L_r$.
Only
$2\,000$
out of
$625^2$
values
($0.512\%$)
need to be computed and~stored. Then,~each value of
$g_N$
is looked up
in~$L_r$.
Here,~as
$y$
and~$w$
are both divisible
by~$5$,
only one value out of
$25$
($4\%$)
needs to be computed and searched~for.
\end{ttt-s}

\begin{ttt-s}
{\bf The Case $5 | x$.}
We call this case ``S'' and use the letter~S at a prominent position in the naming of the relevant files of source~code. S~stands for ``Sonderfall'' which means ``exceptional~case''. It is not known whether there exists a solution belonging to~case~S.

Here, we simply interchange both sides of the~equation. Define
$$f_S (z, w) := z^4 + 4 w^4 \qquad {\rm and} \qquad g_S (x, y) := x^4 + 2 y^4.$$
As
$x$
and
$y$
are divisible by
$5$,
we get
$x^4 + 2 y^4 \equiv 0 \pmod {625}$
and, therefore,
$z^4 + 4 w^4 \equiv 0 \pmod {625}$.

Again,~this congruence allows only
$4 \cdot \varphi(625) = 2\,000$
solutions
$(z, w) \in \left( \bbZ/625\bbZ \right)^2$
and these pairs are easily computable,~too. The~fourth roots
of~$(-4)$
in~$\bbZ/625\bbZ$
are
$\pm181$
and~$\pm183$.
For~each
$x \in \bbZ/625\bbZ^{\textstyle *}$,
one has to consider
$z := (\pm 181x \mod 625)$
and~$z := (\pm 183x \mod 625)$.

We store the values of
$f_S$
into the
set~$L_r$.
Then,~we search
through~$L_r$
for the values
of~$g_S$.
As~above, only
$2\,000$
out of
$625^2$
values need to be computed and~stored and one value out of
$25$
needs to be computed and searched~for.
\end{ttt-s}

\subsection{The prime $2$}

\begin{ttt-s}
Any primitive solution is of the form that
$x$
and~$z$
are odd while
$y$
and~$w$
are~even.
\end{ttt-s}

\begin{ttt-s}
In~case~S, there is no way to do better than that as both
$f_S$
and~$g_S$
represent
$(r \mod 2^k)$
for~$k \geq 4$
if and only
if~$r \equiv 1 \pmod {16}$.

In~case~N, the situation is somewhat~better.
$g_N (y, w) = 4w^4 - 2y^4$~is
always divisible
by~$32$
while
$f_N (x, z) = x^4 - z^4 \equiv 0 \pmod {32}$,
as may be seen by inspecting the fourth roots of unity
modulo~$32$,
implies the
condition~$x \equiv \pm z \pmod 8$.
This~may be used to halve the size
of~$L_r$.
\end{ttt-s}

\subsection{The prime $3$}

\begin{ttt-s}
Looking for further congruence~conditions, a primitive solution must necessarily satisfy, we did not find any reason to distinguish more~cases. But~there are a few more congruences which we used in order to reduce the size of the
sets~$L_r$.

To~explain them, let us first note two theorems on binary quadratic~forms. They~may both be easily deduced from~\cite[Theorems~246 and~247]{Hardy/Wright}.
\end{ttt-s}

\begin{theo-s}
The quadratic forms\/
$q_1 (a,b) := a^2 + b^2$,
$q_2 (a,b) := a^2 - 2 b^2$,
and\/
$q_3 (a,b) := a^2 + 2 b^2$
admit the property~below.

Suppose
$n_0 := q_i (a_0, b_0)$
is divisible by a
prime\/~$p$
which is not represented by
$q_i$.
Then,~$p | a_0$
and~$p | b_0$.
\end{theo-s}

\begin{theo-s}
A prime number\/
$p$
is represented by
$q_1$,
$q_2$,
or~$q_3$,
respectively, if~and only if\/
$(0 \mod p)$
is represented in a non-trivial~way. In~particular,
\begin{iii}
\item
$p$
is represented by\/
$q_1$
if and only if\/
$p = 2$
or
$p \equiv 1 \pmod 4$.
\item
$p$
is represented by\/
$q_2$
if and only if\/
$p = 2$
or
$\big( \frac{2}{p} \big) = 1$.
The~latter means
$p \equiv 1, 7 \pmod 8$.
\item
$p$
is represented by\/
$q_3$
if and only if\/
$p = 2$
or
$\big( \frac{-2}{p} \big) = 1$.
The~latter is equivalent
to\/~$p \equiv 1, 3 \pmod 8$.
\end{iii}
\end{theo-s}

\begin{rem-s}
There is the obvious asymptotic estimate
\begin{align*}
\sharp \{ q_i (a,b) \mid a,b \in \bbZ, q_i (a,b) \in \bbP, q_i (a,b) \leq n \} & \,\sim\, \frac{n}{2 \log n}.  \\
\tag*{Further,} \\
\sharp \{ q_i (a,b) \mid a,b \in \bbZ, |q_i (a,b)| \leq n \} & \,\sim\, C_i \frac{n}{\sqrt{\log n}}
\end{align*}
where
$C_1$,
$C_2$,
and
$C_3$
are constants which can be expressed explicitly by Euler~products.
(For~$q_1$,
this is worked out in~\cite[Satz~(1.8.2)]{Bruedern}. For~the other forms, J.~Br\"udern's argument works in the same way without essential~changes.)
\end{rem-s}

\begin{c81-s}
~\\
In case~N,
$g_N (y, w) = (2w^2)^2 - 2 (y^2)^2 = q_2 (2w^2, y^2)$
where
$q_2$
does not represent the
prime~$3$.
Therefore,~if
$3 | g_N (y, w)$
then
$3 | 2w^2$
and
$3 | y^2$
which implies
$y$
and
$w$
are both divisible
by~$3$.
By~consequence, if
$3 | g_N (y, w)$
then,
automatically,~$81 | g_N (y, w)$.

If
$3 | f_N (x, z)$
but
$81 \notd f_N (x, z)$
then
$f_N (x, z)$
does not need to be stored
into~$L_r$.
Further,~if
$3|x$
and~$3|z$
then 
$f_N (x, z)$
does not need to be stored, either, as it cannot lead to a primitive~solution. This~reduces the size of the
set~$L_r$
by a factor
of~$\smash{\frac19 + 4 \cdot \frac13 (\frac13 - \frac1{81}) = \frac{131}{243} \approx 53.9\%.}$

In case~S, the situation is the other
way~round.~$f_S (z, w) = (z^2)^2 + (2w^2)^2 = q_1 (z^2, 2w^2)$
and
$q_1$
does not represent the
prime~$3$.
Therefore,~if
$3 | f_S (z, w)$
then
$3 | z^2$
and
$3 | 2w^2$
which implies that
$z$
and
$w$
are both divisible
by~$3$
and~$81 | f_S (z, w)$.

We use this in order to reduce the time spent on~reading.
If~$3 | g_S (x, y)$
but
$81 \notd g_S (x, y)$
or if
$3|x$
and~$3|y$
then
$g_S (x, y)$
does not need to be searched~for. Although~modular operations are not at all fast, the reduction of the number of attempts to read
by~$53.9\%$
is highly~noticeable.
\end{c81-s}

\subsection{Some more hypothetical improvements}

\begin{ttt-s} ~
\begin{iii}
\item
In the argument for case~N given above,
$p = 3$
might be replaced by any other prime
$p \equiv 3, 5 \pmod 8$.

In~case~S, the same argument as above works for every
prime~$p \equiv 3 \pmod 8$.
For~primes
$p \equiv 5 \pmod 8$,
the strategy could be~reversed.
$q_3$
is a binary quadratic form which
represents~$(0 \mod p)$
only in the trivial~manner. Therefore,~if
$p | g_S (x, y)$
then
$p | x$
and~$p | y$.
It~is unnecessary to store
$f_S (z, w)$
if
$p|z$
and~$p|w$
or if
$p | f_S (z, w)$
but~$p^4 \notd f_S (z, w)$.
\item[i${}^\prime$) ]
Each argument mentioned may be extended to {\em some\/}
primes~$p \equiv 1 \pmod 8$.
For~example, in case~N, what is actually needed is that
$2$
is not a fourth power
modulo~$p$.
This~is true, e.g., for
$p = 17$,
$41$,
and~$97$,
but not for
$p = 73$
and~$89$.
\item \label{13-17}
$f_N$
and~$f_S$
do not represent the residue classes of
$6$,
$7$,
$10$,
and~$11$
modulo~$17$.
$g_N$~and~$(-g_S)$
do not represent
$1$,
$3$,
and~$9$
modulo~$13$.
This~could be used to reduce the load for writing as well as~reading.
\end{iii}
\end{ttt-s}

\begin{rems-s} ~
\begin{abc}
\item
We did not implement these improvements as it seems the gains would be marginal or the cost of additional computations would even dominate the~effect. It~is, however, foreseeable that these congruences will eventually become valuable when the speed of the CPU's available will continue to grow faster than the speed of~memory. Observe~that alone the congruences noticed in~a) could reduce the amount of data to be stored into
$L$
to a size asymptotically less than
$\varepsilon B^2$
for
any~$\varepsilon > 0$.
\item
For every prime
$p$
different from
$2$,
$5$,
$13$,
and~$17$,
the quartic forms
$f_N$,
$g_N$,
$f_S$,
and~$g_S$
represent all residue classes
modulo~$p$.
This~means, \ref{13-17} may not be carried over to any further~primes.

This can be seen as follows.
Let~$b$
be equal to
$f_N$,
$f_S$,
$g_N$,
or~$g_S$.
$(0 \mod p)$~is
represented
by~$b$,~trivially.
Otherwise,~$b (x, y) = r$
defines an affine curve
$C_r$
of genus three with at most four points on the infinite~line. The~Weil conjectures~\cite[Corollaire~3 du Th\'eor\`eme~13]{Weil} imply that
$[(p + 1 - 6 \sqrt{\mathstrut p}) - 4]$
is a lower bound for the number of
$\bbF_p$-rational
points
on~$C_r$.
This~is a positive number as soon
as~$p \geq 43$.
In~this
case, every residue
class~$(r \mod p)$
is represented, at least,~once.

For~the remaining primes up
to~$p = 41$,
an experiment shows that all residue classes
modulo~$p$
are represented by
$f_N$,
$f_S$,
$g_N$,
as well
as~$g_S$.
\end{abc}
\end{rems-s}

\section{A 64~bit based implementation of the algorithm}

\begin{ttt}
We migrated the implementation of Algorithm~H from a 32~bit processor to a 64~bit~processor. This means, the new hardware supports addition and multiplication of 64~bit~integers. Even~more, every operation on (unsigned) integers is automatically
modulo~$2^{64}$.

From this, various optimizations of the implementation described in~\cite{Joerg_und_Stephan} are almost~compelling. The~basic idea is that 64 bits should be enough to define hash value and control value, two integers significantly less than
$2^{32}$
which should be independent on each other, by selection of bits instead~of using (notoriously slow) modular~operations.

Note, however, that the congruence conditions
modulo~$2$
imposed imply that
$x^4 \equiv z^4 \equiv 1 \pmod {16}$
and
$2 y^4 \equiv 4 w^4 \equiv 0 \pmod {16}$.
This~means, the four least significant bits of
$f$
and~$g$
may not be used as they are always the~same.
\end{ttt}

\begin{ttt}
The description of the algorithm below is based on case~S, case~N being completely~analogous.\smallskip

{\bf Algorithm H64.}
\begin{III}
\item[I. {\em Initialization.} ] \addtocounter{III}{1}
Fix
$B := 10^8$.
Initialize a hash~table of
$2^{27} = 134\,217\,728$
integers, each being 32~bit~long. Fix~the page~prime
$p_p  := 200\,003$.

Further, define two functions, the
{\em hash~function}~$h$
and the
{\em control~function}~$c$,
which map 64 bit integers to 27 bit integers and 31 bit integers, respectively, by selecting certain~bits. Do~not use any of the bits twice to ensure
$h$
and
$c$
are independent on each other and do not use the four least significant~bits.
\item {\em Loop.}
Let
$r$
run from
$0$
to 
$p_p - 1$
and execute steps \ref{Writing} and~\ref{Reading} for
each~$r$.
\begin{ABC}
\item {\em Writing.}
Build up the hash table, which is meant to encode the set
$L_r$,
as~follows.
\begin{abc}
\item \label{zw_auflisten}
Find all pairs
$(z, w)$
of non-negative integers less than or equal to
$B$
which satisfy
$z^4 + 4 w^4 \equiv r \pmod {p_p}$
and all the congruence-conditions for primitive solutions, listed~above. (Make~systematic use of the Chinese remainder~theorem.)
\item
Execute steps \ref{erster_w} and \ref{letzter_w} below for each such pair.
\begin{iii}
\item \label{erster_w}
Evaluate
$f_S (z, w) := (z^4 + 4 w^4 \mod 2^{64})$.
\item \label{letzter_w}
Use the hash~value
$h (f_S (z, w))$
and linear probing to find a free place in the hash~table and store the control~value
$c (f_S (z, w))$~there.
\end{iii}
\end{abc}
\item {\em Reading.}
Search within the hash~table, as~follows.
\begin{abc}
\item \label{xy_auflisten}
Find all pairs
$(x, y)$
of non-negative integers less than or equal to
$B$
which satisfy
$x^4 + 2 y^4 \equiv r \pmod {p_p}$
and all the congruence~conditions for primitive solutions, listed~above. (Make~systematic use of the Chinese remainder-theorem.)
\item
Execute steps \ref{erster_r} and \ref{letzter_r} below for each such~pair.
\begin{iii}
\item \label{erster_r}
Evaluate
$g_S (x, y) := (x^4 + 2 y^4 \mod 2^{64})$
on all points found in step 
\ref{xy_auflisten}.
\item \label{letzter_r}
Search for the control~value
$c (g_S (x, y))$
in the hash table, starting at the hash~value
$h (g_S (x, y))$
and using linear probing, until a free position is~found. Report~all hits and the corresponding values
of~$x$
and~$y$.
\end{iii}
\end{abc}
\end{ABC}
\end{III}
\end{ttt}

\begin{rems}[{\rm Some details of the implementation}{}] ~
\begin{iii}
\item
The fourth powers and fourth roots
modulo~$p_p$
are computed during the initialization part of the program and stored into arrays because arithmetic
modulo~$p_p$
is slower than memory~access.
\item
The control value is limited to 31~bits as it is implemented as a signed~integer. We~use the value
$(-1)$
as a marker for an unoccupied place in the hash~table. 
\item
In contrast to our previous programs~\cite{Joerg_und_Stephan}, we do not precompute large tables of fourth powers
modulo~$2^{64}$
because an access to these tables is slower than the execution of two multiplications in a row (at least on our~computer).
\item
It is the impact of the congruences
modulo~$625$,
$8$,
and~$81$,
described above, that the set of
pairs~$(y, w)$
$[(x, y)]$
to be read is significantly bigger than the set of
pairs~$(x, z)$
$[(z, w)]$
to be~written. They~differ actually by a factor of
$\frac{625^2}{2\,000 \cdot 25}\cdot\frac{243}{112}\cdot2 \approx 33.901$
in case~N and
$\frac{625^2}{2\,000 \cdot 25}\cdot\frac{112}{243} \approx 3.601$
in case~S.

As~a consequence of this, only a small part of the running-time is spent on~writing. The~lion's share is spent on unsuccessful searches
within~$L$.
\end{iii}
\end{rems}

\begin{remso}[{\rm Post-Processing}{}] ~
\begin{iii}
\item
Most of the hits found in the hash~table actually do not correspond to solutions of the Diophantine~equation. Hits~indicate only a similarity of bit-patterns.
Thus,~for each pair
of~$x$
and~$y$
reported, one needs to check whether a suitable pair
of~$z$
and~$w$
does~exist. We~do this by recomputing
$z^4 + 4 w^4$
for all
$z$
and~$w$
which fulfill the given congruence~conditions
modulo~$p_p$
and powers of the small~primes.

Although this method is entirely primitive, only about 3\% of the total running-time is actually spent on post-processing. One~reason for this is that post-processing is not called very often, on average only once on about five~pages. For~those pages, the writing part of the algorithm needs to be~recapitulated. This is, however, not time-critical as only a small part of the running-time is spent on writing,~anyway.
\item
An interesting alternative for post-processing would be to apply the theory of binary quadratic~forms. The~obvious strategy is to factorize 
$x^4 + 2 y^4$
completely into prime~powers and to deduce from the decomposition all
pairs~$(a, b)$
such that
$a^2 + b^2 = x^4 + 2 y^4$.
Then,~one may check whether for one of them both
$a$
and~$\frac{b}{2}$
are perfect~squares.
\end{iii}
\end{remso}

\begin{rem}
The migration to a more bit-based implementation led to an increase of the speed of our programs by a factor of approximately~1.35.
\end{rem}

\section[Adaption to the memory~architecture]{Adaption to the memory~architecture of our computer -- generalities}
\setcounter{ttt-s}{0}

\begin{ttt-s}
The factor of~1.35 is less than what we actually hoped~for. For~that reason, we made various tests in order to find out what the limiting bottleneck of our program~is. It~turned out that the major slowdown is the access of the processor to main~memory. 

Our programs are, in fact, doing only two things, integer~arithmetic and memory~access. The~integer execution~units of modern processors are highly optimized circuits and several of them work in parallel inside one processor. They work a~lot faster than main~memory~does. In~order to reach a further improvement, it will therefore be necessary to take the architecture of memory into closer~consideration.
\end{ttt-s}

\subsection{The memory architecture}

\begin{sit-s}
Computer designers try to bridge the gap between the fast processor and the slow memory by building a memory~hierarchy which consists of several cache~levels.

The~cache is a very small and fast memory inside the~processor. The~first cache~level, called L1~cache, of our processor consists of a data~cache and an instruction~cache. Both~are 64~kByte in~size. The~cache manager stores the most recently used data into the cache in order to make sure a second access to them will be~fast.

If the cache~manager does not find necessary data within the L1~cache then the processor is forced to~wait. In~order to deliver data, the cache~management first checks the L2~cache which is 1024~kByte~large. It~consists of 16384 lines of 64~Byte,~each.
\end{sit-s}

\begin{prog-s}
Our~program fits into the instruction~cache,~completely. Therefore,~no problem should arise from~this.

When we consider the data~cache, however, the situation is entirely~different. The~cache~manager stores the 1024 most recently used memory~lines, each being 64~Byte long, within the L1~data~cache.

This strategy is for sure good for many~applications. It~guarantees main memory may be scanned at a high~speed. On the other hand, for our application, it fails~completely. The~reason is that access to our 500~MByte hash~table is completely~random. An~access directly to the L1~cache happens in by far less than 0.1\% of the~cases. In~all other cases, the processor has to~wait.

Even worse, it is clear that in most cases we do not even access the L2~cache. This~means, the cache~manager needs to access main memory in~order to transfer the corresponding memory~line of 64~Byte into the L1~cache. After~this, the processor may use the~data. In~the case that there is no free line available within the L1~cache, the cache~manager must restore old data back to main memory,~first. This~process takes us 60~nanoseconds, at~least, which seems to be short, but the processor could execute more than 100~integer~instructions during the same~time.

The philosophy for further optimization must, therefore, be to adapt the programs as much as possible to our hardware, first of all to the sizes of the L1 and L2~caches.
\end{prog-s}

\begin{prob-s}
Unfortunately, the whole memory hierarchy is invisible from the point of view of a higher programming~language, such~as~C, since such languages are designed for being machine-independent. Further,~the hardware executes the cache~management in an automatic~manner. This~means, even by programming in assembly, one cannot control the cache completely although some new assembly instructions such~as~{\tt prefetch} allow certain direct~manipulations.
\end{prob-s}

\begin{sol-s}
A practical way, nonetheless to gain some influence on the memory~hierarchy, is to rearrange the algorithm in an apparently nonsensical manner, thereby making memory access less~chaotic. One~may then hope that the automatic management of the cache, when confronted with the modified algorithm, is able to react more~properly. This~should allow the program to run~faster. 
\end{sol-s}

\subsection{Our first trial}

\begin{ttt-s}
Our first idea for this was to work with two arrays instead of~one.\medskip\pagebreak[3]

{\bf Algorithm~M.}
\begin{iii}
\item
Store~the values
of~$f$
into an array and the values
of~$g$
into a another~one. Write~successively calculated values into successive~positions. It is clear that this part of the algorithm is not troublesome as it involves a~linear memory access which is perfectly supported by the memory~management.
\item
Then,~use Quicksort in order to sort both~arrays. In~addition to being fast, Quicksort~is known to have a good memory~locality when large arrays are~sorted.
\item
In a final step, search for matches by going linearly through both arrays as in~Mergesort.
\end{iii}
\end{ttt-s}

\begin{rem-s}
Unfortunately, the idea behind Algorithm~M is too simple to give it any chance of being superior to the previous algorithms. However, it is a worthwhile~experiment. Indeed,~our implementation of Algorithm~M causes at least 30~times more memory transfer compared with the previous programs but, actually, it is only three times~slower. This~indicates that our approach is~reasonable.
\end{rem-s}

\section{Hashing with partial presorting}

\subsection{The algorithm}

\begin{ttt-s}
Our final algorithm is a combination of sorting and hashing. An~important aspect of it is that the sorting~step has to be considerably faster than the Quicksort algorithm. For~that reason, we adopted some ideas from linear-time sorting algorithms such as Radix~Sort or Bucket~Sort.
\end{ttt-s}  

\begin{ttt-s}
The algorithm works as~follows. Again,~the description is based on case~S, case~N being~analogous.\smallskip

{\bf Algorithm H64B.}
\begin{III}
\item[I. {\em Initialization.} ] \addtocounter{III}{1}
Fix
$B := 10^8$.
Initialize a hash~table
$H$
of
$2^{27} = 134\,217\,728$
integers, each being 32~bit~long. Fix~the page~prime
$p_p  := 200\,003$.

In~addition, initialize
$1024$
auxiliary arrays
$A_i$
each of which may contain
$2^{17} = 131\,072$
long~(64~bit)~integers.

Further, define two functions, the
{\em hash~function}~$h$
and the
{\em control~function}~$c$,
which map 64 bit integers to 27 bit integers and 31 bit integers, respectively, by selecting certain~bits. Do~not use any of the bits twice to ensure
$h$
and
$c$
are independent on each other and do not use the four least significant~bits.

Finally,
let~$h^{(10)}$
denote the function mapping 64 bit integers to integers
within~$[0, 1023]$
which is given by the ten most significant bits
of~$h$.
In~other words, for
every~$x$,
$h^{(10)} (x)$ is the same as
$h (x)$
shifted to the right by 17~bits.
\item {\em Outer Loop.} \label{Outer_Loop}
Let
$r$
run from
$0$
to 
$p_p - 1$
and execute \ref{Writing} and~\ref{Reading} for
each~$r$.
\begin{ABC}
\item {\em Writing.} \label{Writing}
Build up the hash table, which is meant to encode the set
$L_r$,
as~follows.
\begin{abc}
\item {\em Preparation.}
Find all pairs
$(z, w)$
of non-negative integers less than or equal to
$B$
which satisfy
$z^4 + 4 w^4 \equiv r \pmod {p_p}$
and all the congruence-conditions for primitive solutions, listed~above. (Make~systematic use of the Chinese remainder~theorem.)
\item {\em Inner Loop.} \label{zw_auflisten_b}
Execute steps \ref{erster_wb} -- \ref{letzter_wb} below for each such~pair.
\begin{iii}
\item \label{erster_wb}
Evaluate
$f_S (z, w) := (z^4 + 4 w^4 \mod 2^{64})$.
\item \label{puffer_erstmalig_linear_fuellen}
Do~not store
$f_S (z, w)$
into the hash~table, immediately.
Put
$i := h^{(10)} (f_S (z, w))$,
first.
\item \label{letzter_wb}
Add
$f_S (z, w)$
to the
auxiliary~array~$A_i$.
Maintain~$A_i$
as an unordered list, i.e.~always write to the lowest unoccupied~address.

If there is no space left
in~$A_i$
then output an error message and abort the~algorithm.
\end{iii}
\item {\em Storing.} \label{hash_blockweise_aufbauen}
Let
$i$
run from
$0$
to~$1023$.
For~each~$i$
let~$j$
run through the addresses occupied
in~$A_i$.

For~fixed
$i$
and~$j$,
extract from the 64~bit integer
$\smash{A_i [j]}$
the 27~bit hash~value
$\smash{h (A_i [j])}$
and the 31~bit
control~value~$\smash{c (A_i [j])}$.

Use~the hash-value
$\smash{h (A_i [j])}$
and linear probing to find a free place in the hash~table and store the control-value
$\smash{c (A_i [j])}$~there.
\item {\em Clearing up.} \label{leeren}
Clear the auxiliary arrays
$A_i$
for all
$\smash{i \in [0, 1023]}$
to make them available for~reuse.
\end{abc}
\item {\em Reading.} \label{Reading}
Search within the hash~table, as~follows.
\begin{abc}
\item {\em Preparation.}
Find all pairs
$(x, y)$
of non-negative integers less than or equal to
$B$
which satisfy
$x^4 + 2 y^4 \equiv r \pmod {p_p}$
and all the congruence~conditions for primitive solutions, listed~above. (Make~systematic use of the Chinese remainder-theorem.)
\item {\em Inner Loop.} \label{xy_auflisten_b}
Execute steps \ref{erster_rb} -- \ref{letzter_rb} below for each such~pair.
\begin{iii}
\item \label{erster_rb}
Evaluate
$g_S (x, y) := (x^4 + 2 y^4 \mod 2^{64})$.
\item \label{lesepuffer_aufbauen}
Do not look up
$g_S (x, y)$
in the hash~table,~immediately. Put
$i := h^{(10)} (g_S (x, y))$,
first.
\item \label{letzter_rb}
Add
$g_S (x, y)$
to the
auxiliary~array~$A_i$.
Maintain~$A_i$
as an unordered list, i.e.~always write to the lowest unoccupied~address.

If there is no space left
in~$A_i$
then call \smash{d[i]} and add 
$g_S (x, y)$
to~$A_i$,
afterwards.
\end{iii}
\item {\em Searching. Clearing all buffers.} \label{alle_lesepuffer_leeren}
Let
$i$
run from
$0$
to~$1023$.
For
each~$i$,
call~\smash{d[i]}.

When this is finished, terminate the~algorithm.\addtocounter{abc}{1}
\item[{\em Subroutine\/} d{[}i{]}) ] {\em Clearing a buffer.} \label{lesepuffer_zwischendurch_in_bloecken_leeren}
Let~$j$
run through the addresses occupied
in~$A_i$.
For~fixed~$j$,
search for the control~value
$\smash{c (A_i [j])}$
within the hash~table
$H$,
starting at the hash~value
$\smash{h (A_i [j])}$
and using linear probing, until a free place is~found. Report~all hits and the corresponding values
of~$x$
and~$y$.

Having~done this, declare
$A_i$
to be~empty.
\end{abc}
\end{ABC}
\end{III}
\end{ttt-s}\pagebreak[3]
 
\begin{rem-s}
The auxiliary~arrays
$A_i$
play the role of a~buffer. Thus,~one could say that we introduced some buffering into the management of the
hash~table~$H$.
However,~this description misses the~point.

What is more important is that the values
of~$f_S$
to be stored into
$L_r$
are partially sorted according to the 10~most significant bits of
$h (f_S (z, w))$
by putting them into the
auxiliary~arrays~$A_i$.
When~the hash~table is then built~up, the records arrive almost in~order. The~same is true for~reading.

What~we actually did is, therefore, to introduce some {\em partial presorting\/} into the management of the hash~table.
\end{rem-s}

\begin{rem-s}
It is our experience that each auxiliary array carries more or less the same~load. In~particular, in step \ref{Outer_Loop}\ref{Writing}b.\ref{letzter_wb}, when the buffers are filled up for writing, a buffer overflow should never~occur. For~this reason, we feel free to treat this possibility as a fatal~error.
\end{rem-s}

\subsection{Running-Time}

\begin{ttt-s}
Algorithm~H64B uses about three times more memory than our previous algorithms but our implementation runs almost three times as~fast. It was this~factor which made it possible to attack the bound
$B = 10^8$
in a reasonable amount of~time.

The~final version of our programs took almost exactly 100~days of CPU time on an AMD~Opteron 248~processor. This~time is composed almost equally of 50~days for case~N and 50~days for case~S. The~main computation was executed in parallel on two machines in February and March,~2005.
\end{ttt-s}

\begin{why-s}
To answer this question, one has to look at the impact of the~cache. For~the old program, the cache~memory was mostly~useless. For~the new program, the situation is completely~different.

When the auxiliary arrays are filled in step \ref{Outer_Loop}\ref{Writing}b.\ref{puffer_erstmalig_linear_fuellen} and \ref{Outer_Loop}\ref{Reading}b.\ref{lesepuffer_aufbauen}, access to these arrays is~linear. There~are only 1024~of~them which is exactly the number of lines in the L1~cache. When~an access does not hit into that innermost cache then the corresponding memory~line is moved to it and the next seven accesses to the same auxiliary array are accesses to that~line. Altogether,~seven of eight memory accesses hit into the L1~cache.

When an auxiliary array is emptied in step \ref{Outer_Loop}\ref{Writing}b.\ref{leeren} or \ref{Outer_Loop}\ref{Reading}b.d[i]), the situation is~similar. There~are a high number of accesses to a very short segment of the hash~table. This~segment fits completely into the L2~cache. It~has to be moved into that cache,~once. Then, it can be used many~times. Again,~access to the auxiliary array is linear and a~hit into the L1~cache occurs in seven of eight~cases.

All in all, for Algorithm~H64B, most memory~accesses are hits into the~cache. This~means, at the cost of some more data transfer altogether, we achieved that main~memory may be mostly used at the speed of the~cache.
\end{why-s}

\frenchspacing\addtolength{\labelsep}{-2.5pt}
\def\leftmark{\scshape References}%

\renewcommand{\thefootnote}{\fnsymbol{footnote}}
\footnotetext[0]{version of June 28${}^{\rm th}$, 2005}
\end{document}